\def\hyperreal{{^*{\real}}}
\def\real{{\tt I\kern-.2em{R}}}
\def\hypernat{{^*{\nat }}}
\def\Hyper#1{\hyper {\eskip #1}}
\def\eskip{\hskip.25em\relax}
\def\nat{{\tt I\kern-.2em{N}}}
\def\hyper#1{\ ^*\kern-.2em{#1}}
\def\Hyper#1{\hyper {\eskip #1}}
\def\power#1{{{\cal P}(#1)}}
\def\qed{{\vrule height6pt width3pt depth2pt}\par\medskip}

\def\b#1{{\bf #1}}
\def\st#1{{\tt st}(#1)}
\def\St#1{{\tt St}(#1)}
\def\real{{\tt I\kern-.2em{R}}}
\def\nat{{\tt I\kern-.2em{N}}}

\def\realp#1{{\tt I\kern-.2em{R}}^#1} 
\def\natp#1{{\tt I\kern-.2em{N}}^#1}
\def\hyper#1{\ ^*\kern-.2em{#1}}

\def\hyperreal{{^*{\real}}}
\def\hyperrealp#1{{\tt ^*{I\kern-.2em{R}}}^#1} 
\def\hypernat{{^*{\nat }}}
\def\hypernatp#1{{{^*{{\tt I\kern-.2em{N}}}}}^#1} 
\def\eskip{\hskip.25em\relax}

\def\Hyper#1{\hyper {\eskip #1}}
\def\leaderfill{\leaders\hbox to 1em{\hss.\hss}\hfill}
\def\srealp#1{{\rm I\kern-.2em{R}}^#1}
\def\sp{\vert\vert\vert} 

\def\power#1{{{\cal P}(#1)}}

\def\qed{{\vrule height6pt width3pt depth2pt}\par\medskip}
\def\pars{\par\smallskip}
\def\parm{\par\medskip}
\def\r#1{{\rm #1}}
\def\b#1{{\bf #1}}
\def\ref#1{$^{#1}$}

\def\m@th{\mathsurround=0pt}
\def\rightarrowfill{$\m@th \mathord- \mkern-6mu \cleaders\hbox{$\mkern-2mu 
\mathord- \mkern-2mu$}\hfil \mkern-6mu \mathord\rightarrow$}
\def\leftarrowfill{$\mathord\leftarrow
\mkern -6mu \m@th \mathord- \mkern-6mu \cleaders\hbox{$\mkern-2mu 
\mathord- \mkern-2mu$}\hfil $}
\def\noarrowfill{$\m@th \mathord- \mkern-6mu \cleaders\hbox{$\mkern-2mu 
\mathord- \mkern-2mu$}\hfil$}
\def\orgate{$\bigcirc \kern-.80em \lor$}
\def\andgate{$\bigcirc \kern-.80em \land$}
\def\inverter{$\bigcirc \kern-.80em \neg$}

\magnification=\magstep1 
\tolerance 10000
\baselineskip  14pt
\hoffset=.25in
\hsize 6.00 true in
\vsize 8.85 true in
\centerline{{\bf The GGU-model and Generation of Developmental Paradigms}}\par\medskip
\centerline{Robert A. Herrmann}
\centerline{4 MAY 2006, 16 FEB 2014}
\bigskip
\noindent {\bf 1. Introduction.}\par\medskip
The most basic aspects of the General Grand Unification Model (GGU-model) are delineated in [6].  In this paper, a definition for a refined development paradigm is  formalized and a new method to obtain correspond ultrawords is established. The notion of the ``logic-system signature'' for a logic-system generated by a scientific theory consequence operator ${\rm S_N}$ (in any of its forms) [4] is detailed. It is shown that the extended standard part operator, $'\tt St$, is a finite consequence operator. \par\medskip 

\noindent {\bf 2. Informal Developmental Paradigms}. \par\medskip 
A nonempty ``alphabet'' ${\cal A}$ is considered as a finite or denumerable set. The informal language L generated from $\cal A$ [7] has the property that $\vert {\rm L} \vert = \vert \bigcup \{{\cal A}^n\mid (n \in \nat)\land (n > 0)\}\vert.$  The denumerable language L can be considered as informally presented and then the following defined sequences embedded as sequences into the G-structure [7], or you can begin with members of ${\cal E} = \b L$ and construct these developmental paradigms. These developmental paradigms yield the quasi-physical event sequences. The former method is used in what follows. Further, for simplicity, consider a ``beginning'' frozen segment $\rm F$ [7]. These frozen segments correspond to the notion of a ``frozen-frame'' in [3]. \par\smallskip 

The idea is to employ the notion of ``finite'' choice to characterize, at least, partially the developmental paradigms for nonempty countable $\rm D\subset L.$ Let $\nat$ denote the natural numbers. For each $\rm a,\ b \in \nat,\ a \leq b,\ [a,b]=\{x \mid (a \leq x\leq b)\land(x \in \nat)\}$ and, as usual, symbol $\rm f|A$ denotes the restriction of the function f to a subset A of its domain. As usual, $\rm D^{[a,b]}$ denotes the set of all functions on $\rm [a,b]$ into $\rm D.$ \par\medskip

\noindent{\bf Definition 2.1.} (1) For $1 \in \nat,$ let $\rm {\cal X}_1 = \{f \in D^{[0,1]}\mid (f(0) = F)\land (f(1) \in D)\}.$\par\smallskip
(2) Assume that $\rm {\cal X}_n$ has been defined. Let $\rm {\cal X}_{n+1} = \{f\in D^{[0,n+1]}\mid (f|[0,n] \in {\cal X}_n)\land (f(n+1)  \in D)\}.$ \par\smallskip
(3) Define $\rm {\cal X} =\bigcup\{{\cal X}_n\mid n \in \nat\}.$\par\smallskip
(4) The actual developmental paradigms $\rm DP$ for a particular $\rm F$ are $\rm DP = \{f\in D^A\mid (A = \nat)\land (f(0)=F)\land (\forall n ((n>0)\land (n \in A) \to  (f|[0,n] \in {\cal X}\}$. Note: For $\rm A = \nat,\ n >0$, if $\rm f\in D^A,$ then $\rm f|[0,n] \in {\cal X}$ if and only if $\rm f|[0,n]\in {\cal X}_n.$\par\medskip

All of this is extended to the hyperfinite when embedded into the nonstandard structure.
In [3], a ``master'' event sequence is used in an attempt to model Definition 2.1 in a reasonably comprehensible manner using various constructive illustrations. Notice that such a master event sequence is a member of $\rm DP$ under Definition 2.1. Obviously, Definition 2.1 is not the only way to obtain developmental paradigms. Indeed, a simple induction proof shows that $\rm \{f \in D^A\mid (A = \nat)\land(f(0) = F)\}= DP.$ Clearly, $\rm DP \subset \{f \in D^A\mid (A = \nat)\land(f(0) = F)\}.$ Let $\rm d \in \{f \in D^A\mid (A = \nat)\land(f(0) = F)\}.$ Then $\rm d(0) = F,$ and $\rm d(1) \in D$ imply $\rm d|[0,1]\in D^{[0,1]}$ and $\rm d|[0,1]\in {\cal X}_1.$ Suppose that $\rm d|[0,n] \in {\cal X}_n.$ Then $\rm d|[0,n+1]$ means that $\rm d|[0,n]$ and $\rm d|{n+1}$ and $\rm d(n+1) \in D.$ Thus, since $\rm d|[0,n+1] \in D^{[0,n+1]},\ d(0) = F,\ d|[0,n] \in {\cal X}_n,$ it follows that $\rm d|[0,n+1] \in {\cal X}_{n+1}.$ Therefore, by induction for each $\rm n>0,\ n \in \nat,\  d|[0,n] \in {\cal X}_n,$ and $\rm d(0) = F$. Hence, $\rm d \in DP.$ For a specific D, among the members of DP are the members used for the GGU-model. \pars 

If Definition 2.1 is restricted to the ``potential'' infinite, then step (3) and (4) are not included. It is step (4) that some might consider as requiring the Axiom of Choice in that members of each $\rm {\cal X}_n$ that are not completely specified are employed in the set-theoretic definition. This somewhat constructive way to define the DP is used to indicate that, at the least, major portions of the basic definition can be obtained via finite choice. This type of finite characterization, when extended to the NSP world, allows for further interesting observations. For example, the event hyperfinite sequences can have ultranatural events not merely associated with nonstandard primitive time but also can have them at standard moments of primitive time. This additional property has not been discussed in [3]. Moreover, this possibility also leads to ultrawords and ultimate ultrawords to which ultralogics can be applied since such developmental paradigms can be considered as of the $d'$ type discussed in section 9.1 in [7]. \par\medskip

\noindent{\bf 3. Refined Developmental Paradigms.}\par\medskip

The following conventions are used. With certain exceptions, each member of the informal set-theory employed is represented by roman fonts.  In most cases, natural numbers, integers and rational numbers and sequences of these are represented by math-italics in both the informal and formal structures. The illustrations in [3] for generating event sequences using Definition 2.1 do not correspond to the actual technical definition that appears in Chapter 7 in [7] except under a specific restriction. (Note: The usual structure now employed for what follows is the Extended Grundlegend Structure (EGS) ${\cal Y}_1$ as defined in [7, p. 70]. The ground set for the standard superstructure is the set of atoms 
$A_1 \cup A,\ A_1 \cap A =\emptyset$, where $A_1$ is isomorphic to the natural numbers and $A$ is isomorphic to the real numbers $\real,$ and, hence, the set $A$ is usually denoted by $\real$.) Definition 2.1 and these illustrations if restriction to a small primitive time interval $[a,b)$ do correspond to those used in [7]. In that case, the actual complete developmental paradigm would be a countable collection of such developmental paradigms for each $[a,b)$. This would technically require that, for applications such as discussed in [3], an additional collection of ultimate ultrawords be considered via Theorem 7.3.4 in [7]. However, the complete developmental paradigm as the denumerable union of denumerably many sets can also be considered as a denumerable sequence in primitive ``time'' when the Axiom of Choice is assumed. In this case, the complete developmental paradigm can be generated by a basic ultraword and the ultralogic $\Hyper {\bf S}.$ \pars

The method devised in Chapter 7 of [7] to analyze a developmental paradigm is significant and should be used since it yields the greatest control and, in the EGS, displays the ultranatural events. Requiring that the denumerable union of denumerably many objects be denumerable is not necessary if the notion of the developmental paradigm is simply defined via different denumerable sets of primitive identifiers. These notions are now formalized within the standard EGS.\par\smallskip
It can be assumed that what follows is the result of an embedding into the standard superstructure of the informal objects. Let $\bf Z$ denote the integers and consider ${\bf Z}\times \nat.$ \par\medskip

\noindent{\bf Definition 3.1.} For each $(i,j),\ (p,m) \in {\bf Z} \times \nat$, where 
$(a,b) = \{\{a\},\{a,b\}\},$ let $\preceq$ be defined as follows:\par\smallskip
(1) if $i < p,$ then $(i,j) \prec (p,m).$\par\smallskip
(2) If $i = p$ and $j <m,$ then $(i,j) \prec (p,m),$ and\par\smallskip
(3) $(i,j) = (p,m)$ if and only if $i = p$ and $j=m.$\par\smallskip 
\noindent The binary relation $\preceq$ yields a simple order in ${\bf Z} \times \nat.$ \par\medskip
The intervals $[a,b)$ employed in [7, p. 61] may be replaced with the following specifically defined intervals. Consider nonzero $K \in \nat.$ When $(i,j)$ appears as a subscript, it is often written as $ij$. For each $i \in {\bf Z},$ let $c_{i} =  i/K.$ For the rational numbers \b Q and each $i \in \b Z$, let $[c_i,c_{i+1}) = \{x\mid(x \in \b Q)\land (c_i \leq x < c_{i + 1})\}.$ For such $i \in\b Z$, partition $[c_i,c_{i+1})$, in the same manner as done in [7, p. 61], by a denumerable increasing sequence of partition points $t_{ij}$, $j \in \nat,$ such that $t_{i0} = c_i,$ $t_{ij} \in [c_{i}, c_{i+1})$ and $\lim_{j \to \infty} t_{ij} \to c_{i+1}.$ For various $i \in \b Z$, the set of rational numbers $\{t_{ij}\mid j \in \nat\}$ models a ``primitive (time) interval.'' For example, let $t_{ij} = (1/K)(i + 1 - 1/2^j).$ For applications, one might employ a finite sequence $\{[0,c_1)\ldots, [c_j,c_{j+1})\}$ of such intervals or partition $[0, +\infty)$, $(-\infty,0),$ or $(-\infty, +\infty)$  using collections of such intervals. If the collection of primitive intervals is nonempty and finite, then there are denumerably many partition points. If the collection of primitive intervals is infinite, then, using the Axiom of Choice, there are denumerably many partition points. If $t_{ij}, \ t_{pm}$ are any of these constructed partition points, then  $t_{ij} \leq t_{pm}$ (the standard rational number simple order) if and only if $(i,j) \preceq (p,m).$ \par\smallskip

By construction, $r \in \b Q$ is a partition point if and only if $r$ corresponds to a frozen segment $\b F_r.$ A major aspect associated with applications is the difference between primitive and observer time.  All of the results in [7] that deal with developmental paradigms are relative to countably many collections of frozen segments. Although for certain applications the actual physical events may be repeated relative to primitive time, the construction of the developmental paradigm allows $\nat$ to be mapped bijectively onto $[a,b).$  Only intervals of the form $[a,b)$ are considered in [7]. Then,  for a countable collection of such partitioned intervals, there is an ultimate ultraword that generates, for each $[a,b)$ interval, the  appropriate ultraword from which each interval's developmental paradigm is obtained.  \par\smallskip

For this refined approach, the use of the $\nat$ notation can be retained under the view that there is a bijection from the set of all partition points onto $\nat$. However, it is a rather trivial matter to re-express each developmental paradigm and its standard frozen segments in terms of the appropriate denumerable subsets of ${\bf Z} \times \nat$ that correspond to the partition points $t_{ij}.$ If this correspondence is employed, then each frozen segment $\b F_{ij}$ corresponds to the partition point $t_{ij}.$ The order $\leq_D$ defined on a developmental paradigm $\b D$ is the simple order induced by $\preceq$ when the developmental paradigm is properly defined. The ``equality'' $=_D$ is set equality. The use of this refined partition point notion yields certain more detailed characteristics for event sequence behavior for the General Grand Unification model (GGU-model) [3, 6].  \par\medskip

\noindent{\bf Examples 3.2.} As an example of the refined use, as in [7], let only members of $[c_i,c_{i+1})$ be considered. In general, let D denote the range for a developmental paradigm for a specific primitive time interval.  [Note: In [7], the range of a general development paradigm range is denoted by d.] Then an embedded standard developmental paradigm is sequentially presented by considering a defining  bijection $f\colon (\{i\}\times \nat) \to \b D.$ Bijections such as this model how members of $\b D$ may be ``grouped together.'' (Note: D now denotes the range.) Each member $\rm d(k) = F$ of D contains a symbol that corresponds to the k (since $\b Q$ is denumerable) and further corresponds under $f$ to a specific moment $(i,j)$ in primitive time. This is a refinement of the construction in [7, p. 61]. Hence, for any $j \in \nat$, and any $p \in \nat,$ such that $j < p, \ f(i,j)<_D f(i,p)$ and $f(i,p) \not= f(i,j).$ Using EGS and *-transfer, this yields that for $\nu \in \hypernat - \nat$ and for $\hyper f$ and each $n \in \nat,\ \hyper {f}(i,n)= f(i,n) \not= \hyper {f}(i,\nu),$ $f(i,n) \hyper {<}_D \hyper{f}(i,\nu),$ and $\hyper {f}(i,\nu) \in \Hyper {\b D}- {\b D}.$ In this application, such objects as $\hyper {f}(i,\nu)$ are ultranatural events. \par\smallskip

Now consider a partition of $(-\infty,0)$ and assume that the developmental paradigm $\b D$ is determined by a bijection $f\colon ({\b Z}^{<0} \times \nat) \to \b D.$ Let $i\in {\b Z}^{<0}.$ Then for any $j \in {\b Z}^{<0}$ such that $i < j,\ f(i,n) <_D f(j,n),\ f(i,n) \not= f(j,n)$ for each $n \in \nat.$  By *-transfer, let $i =\nu \in \Hyper {\b Z}^{<0} - {\b Z}^{<0}.$ It follows that $\hyper {f}(\nu, n) \hyper {<}_D \hyper {f}(j,n),\ \hyper {f}(\nu, n) \not= \hyper {f}(j,n)$ for each $j\in {\b Z}^{<0}$ and $n \in\hypernat.$ Hence, for each $n \in \hypernat,$ $\hyper {f}(\nu, n) \in \Hyper {\b D} - \b D$ since $\nat \subset \hypernat.$ Each of the $\hyper {f}(\nu,n)$ members of $\Hyper {\b D} -\b D$ is called an {\it initial member}. Let $\lambda \in \hypernat - \nat.$ Note that $\b D \subset D_2=\{\hyper f(x,y)\mid (x \in \Hyper {\b Z}^{<0})\land (\nu \leq x) \land (y \in \hypernat)\land (y \leq \lambda)\}\subset \Hyper {\b D}$ and $D_2$ is hyperfinite. For certain applications, one can consider various initial members as corresponding to the same ultranatural event. For the case of the partitioning of $[0,+\infty),$ the same analysis yields members of $\Hyper {\b D}- \b D$ that are termed the {\it final members} in $\Hyper {\b D}.$ For this and other intervals, various members of $D_2$ can be specialized in the sense that they are only members of $D_2$ because of the identifying *-primitive time identifiers. The *-general description $\Hyper f(x,y) \in D_2$ can otherwise be composed of members of $\b L.$ \par\medskip 

In section 10.2 of [7], the hyperfinite choice operator is discussed. This is easily restricted to the ``hyperfinite ordered choice'' operator for hyperfinite subsets of $\Hyper {\b D}.$ By considering any nonempty finite $F\subset \b D$, induction shows there is a $\leq_D$ largest member in $F$. Using this fact, induction shows that for any such finite subset $ F$ there exists a finite ordered choice operator, $O,$ that can be considered as arranging members of $F$ in the proper $\leq_D$ order.  The result of applying $\Hyper {\bf S}$ to an appropriate ultraword yields a hyperfinite set $D'_1$ that contains the embedded developmental paradigm $\bf D$ for interval $[a,b)$ or for intervals such as $(-\infty,0),\ (-\infty,+\infty),\ [0,+\infty)$. Theorems such as  10.1.1 in [7]  are independent from the actual type of partitioning used. They only employ the fact that $\b D$ is denumerable. Using the partition ordering on $\bf D$, where members of $\bf D$ are considered as members of the any of the four primitive interval, the ordering can be restored in the most direct manner by considering the hyperfinite ordered choice operator, $\Hyper {\b O}$. Under this view, the operator $\Hyper {\bf S}$ is composed with $\Hyper {\b O}$ (i.e. $\Hyper {\b O}(\Hyper {\b S}(\{w\})),$ where $w$ is an appropriate ultraword) and this yields all the members of $D_2$ and, hence, $\bf D$ in assigned order.\par    
\medskip 
\noindent{\bf Example 3.3.} Although for this refined definition and its applications a superstructure constructed using the set $\b Q$ as atoms can be used, the superstructure using the reals $\real$ as atoms [7, p. 61] is still employed due to other types of applications. Using the finite power set operator, a finite nonempty collection of finite subsets of $\real$ is a member of $\cal R$ as is the set theoretic union of such a collection. This yields that the union of a nonempty hyperfinite collection of hyperfinite subsets of $\hyperreal$ is a member of $\hyper {\cal R}.$
Consider the standard set $C =\bigcup\{[c_i,c_{i+1})\mid i \in \b Z\},$ where each $t_{ij},\ j \in \nat,$ is as defined previously. The informal set $\rm D$ is denumerable. For the embedded $\b D$, let $f\colon C \to \b D$ be a defining bijection. Using $\rm D$, the informal set $M$ is constructed as done is in [7, p. 66] and embedded into the superstructure. From the definition of $\rm S$ [7, p. 65], it follows that for the set theoretic union $\rm D_1$ having two or more members of a nonempty finite collection of finite subsets of D, there is a word $\rm w \in  M - D$ and $\r D_1 \subset \b S(\{\r w\}).$ Hence, by *-transfer, for any hyperfinite union $F,$ with two or more members, of hyperfinitely many hyperfinite subsets of $\hyper (f[C])= \hyper f[\hyper C]= \Hyper {\b D}$, there exists an ultraword $w' \in \Hyper {\b M} - \Hyper {\b D}$ such that $F \subset \Hyper {\b S}(\{w'\}).$ The members of $\Hyper {\b S}(\{w'\})$ have the same properties described in Theorem 10.1.1 [7, p. 89]. [The method used here to obtain $w'$ for this specific type of developmental paradigm is distinct from that used in [7, p.67].]\par\smallskip

In particular, consider $\lambda < 0 <  \gamma, \ \lambda, \gamma  \in \Hyper {\b Z} -\b Z$ and let $F$ be the (internal) hyperfinite union of the hyperfinite collection of hyperfinite sets $[c_k,c_{k+1}),\ \lambda \leq k \leq \gamma,$ where each $\hyper t_{ij}$ has the property that $0\leq j \leq n_i \in \hypernat - \nat.$ Then $\hyper {f}[F]$ exists, is hyperfinite and $\b D \subset \hyper {f}[F] \subset \Hyper {\b D}.$ Thus, there exists an ultraword $W \in \Hyper {\b M} - \Hyper {\b D}$ such that $\b D \subset \hyper {f}[F] \subset \Hyper {\b S}(\{W\}).$ Applying $\Hyper {\b O}$ to (internal) hyperfinite $\hyper {f}[F]$ yields the ordered $\b D \subset \hyper {f}[F]$. \par\medskip

\noindent{\bf 4. Formalizing Examples 3.2 and 3.3.}\par\medskip

 Recall that the informal general language L is denoted by $\cal W$ in  [7, p. 7]. (The term {\it informal} signifies that general ZFC-set-theory is being used and objects are not considered as members of a superstructure [7].) The set of individuals (atoms, ground set) for the standard superstructure is the set of real numbers $\real$. (This could be changed to the set of rational numbers. However, usually, for the GGU-model, $\real$ is an appropriate set of individuals [7, p. 70.) \pars

For the remainder of this article, the usual notation for developmental paradigms as slightly modified is used. For example, $\r d_q$ denotes the range for a developmental paradigm. Objects in the standard superstructure model ${\cal M} = \langle{\cal R}, \in, =\rangle$ are considered as isomorphically embedded into the superstructure ${\cal Y}_1 = \langle Y, \in,=\rangle.$  Objects in the nonstandard model $\hyper {\cal M} = \langle \Hyper {\cal R}, \in, =\rangle$ are also members of superstructure ${\cal Y}_1$ [7. pp. 22-23].  Since \r L is denumerable, there is a bijection $\r i \colon \r L \to \nat$, where $\nat$ is the set of natural numbers. For $\emptyset \not= H \subset \r i[\r L],\ n \in \nat$, $H^n = H^{[0,n]}$ denotes a set of all functions on $[0,n]$ into $ H$ termed here as {\it partial sequences}.  Let $P = \bigcup \{ (\r i[\r L])^n \mid n \in \nat\}.$ For any $\rm w \in L,$ an equivalence relation is defined relative to how $\r w$ can be formed via the juxtaposition (join) operator [7, pp. 10 - 11.] This equivalence relation separates $P$ into a collection $\cal E$ of equivalences classes. There is a bijection $\theta \colon \r i[{\r L}] \to {\cal E}$, where $\theta(i(\r w ))= [f_0].$ In the usual manner, both $i$ and $\theta$ are extended to set and relation maps. \pars

All members of informal set-theory with constant names that correspond to members produced by the composition $\Theta = (\theta \circ \r i)$ are denoted by the same constant but are written in bold font. All other members of ${\cal Y}_1$ are displayed using math-italics. For example, if $\r C$ is a consequence operator defined on $\r L$, then $\Theta(C) = \b C$ is a consequence operator defined on $\cal E$.  The natural numbers, integers and rational numbers and sets of these entities, except for special notation, for the informal model and standard model are identified and represented by math-italics. \pars

Rather than work in the informal model and map all of the material to $\cal M$ via $\r i$ or $\Theta$, what follows is usually discussed  in terms of informal set-theory or the standard model $\cal M$, where the embedding is not stated specifically. GGU-model primitive time is merely defined as a sequence of rational numbers.  \pars

[For all GGU-model applications as originally presented in [7, Theorem 7.3.1], the developmental paradigm determining functions f and t, as discussed below, are defined on $\r Z \times \nat$ and then the $q$ notion, where $q = 1,\ 2,\ 3$ indicates a restriction of these functions to $\r Z_q \times \nat.$ For $q = 4,$ the indicated functions are the original unrestricted ones. For the t function, the image is $R \subset \b Q.$ Then $\r t^q$ is the appropriate restriction. Hence, the F domain is $R$ and maps $R$ into the informally denoted language L.]\pars

For four sets of integers $\r Z_q,\ q = 1,2,3,4$, the set of rational numbers $\b Q$ and for each $i \in \b Z_q,$ the map $ \r t^q \colon \b Z_q \times \nat \to \b Q$ takes each $j \in \nat$ and yields a basic collection of rational numbers contained in $[c_i, c_{i + 1})= \{ x\mid (c_i \leq x < c_{i+1})\land (x \in \b Q)\},$ where $i \in \b Z$ and $c_{i+1} \in \b Q,$  and similar elementary intervals. Let $R_q = \{\r t^q(i,j)\mid(i \in \b Z_q)\land (j\in \nat)\}\subset \b Q.$  Each member X of a developmental paradigm $\r d_q$ contains a unique rational number identifier taken from $R_q$. This identifier is a member of the word X. Hence, there is a bijection $\r F^q\colon R_q \to \r d_q.$ The values of this bijection, in other articles, are often denoted by $\r F^q_r$ or $\r F^q_{\r t^q(i,j)}$. For the collections of integers, $\r Z_q,\ q=1,2,3,4,$ and for the GGU-model, it is convenient to consider basic intervals as partitioning  the four rational number intervals $q =1 = [0,b],\ (0< b),$ where $[c_m, c_{m+1})$ and $[c_{m-1}, b]$ are employed, and $q = 2 = [0, + \infty),\ q = 3 = (-\infty, 0]$  and $q = 4 = (-\infty, +\infty).$ (Note: The intervals could be real number intervals. However, it appears to be sufficient to consider but the rational numbers.)\pars

 Recall that, for various  $i \in  \r Z_q,$ $j$ varies over the entire set $\nat$ and yields a strictly increasing sequence of rational numbers $\r t^q(i,j) \in [c_i,c_{i+1})$  (and in similar elementary intervals) such that $\lim_{j \to \infty} \r t^q (i,j) = c_{i +1}=\r t^q(i+1,0).$ For each $i \in \r Z_q,$ let each $\r t^{q}(i,j)$ in $R_q$ correspond to  $\r F^q({\r t^{q}({i,j})}) \in \r L$ and the collection of all such $\r F^q({\r t^{q}({i,j})})$ is a developmental paradigm or, when applied to physical-systems, the collection forms an {event sequence}. Informal development paradigms $\r d_q$  correspond to $\b d_q$ in the standard model. Each $\r F^q(\r t^q(i,j))$ is distinct in, at least, one identifying feature -  primitive time. An obvious composition yields a bijection $\r f^{q} = \r F^q \circ \r t^q\colon \b Z_q \times \nat \to \r d_q,  q  =1,2,3,4.$ \pars

(Unless otherwise stated, in all that follows in this article $q =1,2,3,4$, respectively.) In [7, p. 42], it is mentioned that $\r W_{\r j}$ is used to generate specific members of each developmental paradigm d and $\r j,$ as part of the word, is an identifier and can be altered. This should be done for each of the refined developmental paradigms defined on the $q$-intervals. By tick notation, there are symbols in our language for the rational numbers and these correspond to the abstracted notion in the standard model. Hence, each $\b W^q_{t^{q}(i,j)}= \b W^q_r$ is the embedded  statement: ${\rm This\sp description\sp is\sp named\sp}\lceil \r t^{q}{(i,j)}\rceil\sp .$ Descriptions are members of a general language L that not only contains words, in the usual sense, but abstractions of the notion of images and digitized forms of human sensory information. For $(n,m) \not= (r,s)$, the actual event to which  ${\rm xW}^q_{\r t^{q}(m,n)}$  corresponds can be an identical event to which ${\rm xW}^q_{\r t^{q}(r,s)}$  corresponds. The range of $\r t^{q}$ is a denumerable set of rational numbers $R_q$. \pars

The basic construction uses the lexicographic simple order $\preceq$ as defined  on $\b Z_q \times \nat$ and yields an order preserving injection  into the set of all rational numbers $\b Q$, where $\b Q$ carries its standard simple order $\leq$. Trivially, there is an order preserving bijection from $\b Z_q \times \nat$ onto $ R_q$. For an informal developmental paradigm $\r d_q,$ there exists a set ${\cal F}'(\r d_q)$ of all the finite subsets of $\r d_q$ of two or more members. For each $\r G_q \in {\cal F}'(\r d_q),$ first list the members of finite $\r G_q$ from left-to-right without requiring any specific order.  For the next construction, a formal language that is isomorphic to the informal language is employed. \pars

Each $\land$ [resp. $\r F^q({\r t^q(i,j)})$] corresponds to a specific $\rm \sp and\sp$ [resp. a propositional atom that corresponds to a specific word] when  embedded.  This eliminates confusion when $\rm \sp and\sp$ appears in an $\r F^q({\r t^{q}(i,j)})$. Each word is re-expressed by placing $ \land$ between each pair that appear in this first word-form. For example, in terms of propositional atoms if $\r G_q = \{\r F^q({10/3}),\r F^q({7/8}),\r F^q({100/23})\},$ then such a word formed by this construction is $\r w = \r F^q({10/3})\land \r F^q({7/8})\land \r F^q({100/23}) \in \r L$.  Constructions of this form are consistent with the methods used in informal word theory [10, p. 231]. For each $\r G_q$, let $\r w_q$ be one of these constructions and let the injection ${\cal G}\colon {\cal F}'(\r d_q) \to \r L$ map each $\r G_q \in {\cal F}'(\r d_q)$ to $\r w_q.$ Let denumerable set $\r D_q = \{\r w_q \mid({\cal G}(\r G_q) = \r w_q)\land (\r G_q \in {\cal F}'(\r d_q))\}.$ \pars

For the logic-system $S$ [7, pp. 65-66], each $\land $ is interpreted as the $\rm \sp and\sp$ and members of $\r d_q$ are propositions. This corresponds to the above construction. For any $\r w_q = \r F^q(r) \land \cdots \land \r F^q(s) \in \r D_q$,  and the corresponding finite consequence operator, S,  $\{\r F^q(r), \ldots, \r F^q(s)\} \subset \r S(\{\r w_q\}).$ In what follows, although it may not be stated formally, each of the of the standard elements, sets, and relations contained in a formal first-order statement are members of a specific transitive superstructure set $X_p$. Hence, all *-transferred elements, sets and relations are members of  $\hyper X_p.$ There are various ways to obtain ``ultrawords'' [7].  The following method does not specifically use a concurrent relation to obtain ultrawords. Throughout the following, the composition $\b f^{q} = \b F^q \circ t^q$ is employed. \parm     

\noindent{\bf Theorem 4.1.} {\it Consider primitive time interval $1 = [0,b],  b>0.$ It can always be assumed that interval 1 is partitioned into two or more intervals $[c_0,c_1),\ldots$ $ [c_{m-1}, c_m], \ c_m = b, \ m >1,\ m \in \b Z.$  Let ${\b d_1}$ be a developmental paradigm defined on  $R_1 \subset [ 0,b]$. For each infinite $\lambda \in \nat_\infty = \hypernat - \nat,$ there is a $w_\lambda \in \Hyper {\b D}_1$ and an hyperfinite $d^1_\lambda \subset \hyper {\b d}_1$ such that $\b d_1 \subset d^1_\lambda \subset \Hyper {\b S}(\{w^1_\lambda\}).$ Further, $x\in d^1_\lambda$ if and only if there exist an $i\in \Hyper {\b Z}^{\geq 0}, \ 0\leq i\leq m,$ and $j \in \hypernat, \ 0\leq j \leq \lambda,$ such that $x =\Hyper {\b f^1}(i,j).$} \pars

Proof.  Since for each $i,\ 0\leq i < m,\ i \in \b Z,$ $m \in \b Z^{\geq 0},$ the non-negative integers, $\lim_{j \to \infty}t^1(i,j) = c_{i+1}= t^1(i+1,0),$ then $\{\b f^1(i,j)\mid (i \in \b Z)\land(\ 0\leq i < m)\land(j \in \nat)\}\cup \{\b f^1(m,0)\} = \b d_1$.  For each $n \in \nat$, let $H(m,n ) = \{\b f^1(i,j)\mid ( 0\leq i< m)\land (0\leq j \leq n)\land (i \in \b Z)\land(j \in \nat)\} \cup\{\b f^1(m,0)\}.$ Then $H( m,n ) = d^1_n \in  {\cal F}'(\b d_1).$ Further, there is a $\b w^1 \in \b D_1,$ denoted by $\b w^1_n,$ such that $H( m,n ) \subset \b S(\{\b w^1_n\}).$ Consequently, for fixed nonnegative $m \in \b Z,$\par 

$$\forall y((y \in \nat)\to \exists x((x \in \b D_1)\land(H( m,y ) \in  {\cal F}'(\b d_1))\land$$ 
$$(H( m,y ) \subset \b S(\{x\})) )),\eqno (1)$$
holds in $\cal M$ and, hence, the *-transfer of (1) holds in $\hyper {\cal M}.$ This yields\par

$$\forall y((y \in \hypernat)\to \exists x((x \in \Hyper {\b D_1})\land(\hyper H( m,y ) \in \hyper {\cal F}'(\hyper {\b d_1}))\land$$ 
$$(\hyper H( m,y ) \subset \Hyper {\b S}(\{x\})) )),\eqno (2)$$

 Since the union of two hyperfinite sets is hyperfinite, then, for any $y \in \hypernat,$ there is a $ w^1_y \in \Hyper {\b D_1}$ such that  $\hyper H( m,y ) = \{\Hyper {\b f}^1(i,j)\mid ( 0\leq i< m)\land (0\leq j \leq y)\land (i \in \Hyper {\b Z})\land(j \in \hypernat)\}\cup\{\Hyper {\b f}^1(m,0)\}= d^1_y \in \hyper {\cal F}'(\hyper {\b d_1})$, $d^1_y \subset \hyper {\b d_1},$ and $d^1_y \subset \Hyper {\b S}(\{w^1_y\}).$ Notice that if $i \in \b Z$ and $0\leq i\leq m,$ then  $i \in \Hyper {\b Z}$. Also each member of $\Hyper {\cal F}'(\hyper {\b d_1})$ is hyperfinite. \pars

Consider any $\lambda \in \nat_\infty.$  Then the above *-transferred statement holds for $y = \lambda.$ However, if $n \in \nat$, then $n < \lambda$  implies $\Hyper {\b f}^1(i,n) = \b f^1(i,n), \ 0\leq i< m$ and $\Hyper {\b f^1}(m,0) = \b f^1(m,0).$ Hence, $\b d_1 \subset d^1_\lambda.$ From the definition of $H( m,n )$, $x\in d^1_\lambda$ if and only if there exist an $i,\  0\leq i\leq m,$ and $j \in \hypernat, \ 0\leq j \leq \lambda,$ such that $x =\Hyper {\b f^1}(i,j).$ This completes  
the proof. \qed \parm

In the next results, for $q =2,3,4,$ ordered subscript notation is employed and, although the proofs are but modifications of that for Theorem 4.1, separate theorems are presented. \parm
\noindent{\bf Theorem 4.2.} {\it Consider primitive time interval $2 = [0,+\infty).$ Interval $2$ is partitioned into  intervals $[c_i,c_{i+1}),\ i \in \b Z^{\geq 0}.$ Let $\b d_2$ be a developmental paradigm defined on $R_2$. For each infinite $\lambda \in \nat_\infty$ and $\nu \in \b Z_\infty^{\geq 0} = \Hyper {\b Z}^{\geq 0} - \b Z^{\geq 0},$ there is a $w^2_{\nu\lambda} \in \Hyper {\b D}_2$ and an hyperfinite $d^2_{\nu \lambda} \subset \hyper {\b d}_2$ such that
 $\b d_2 \subset d^2_{\nu\lambda}\subset \Hyper {\b S}(\{w^2_{\nu \lambda}\}).$ Further, $x\in d^2_{\nu\lambda}$ if and only if there exist an $i \in \Hyper {\b Z}^{\geq 0},\  0\leq i\leq \nu,$ and $j \in \hypernat, \ 0\leq j \leq \lambda,$ such that $x =\Hyper {\b f}^2(i,j).$ } \pars

Proof. The convergence requirement implies that $\{\b f^2(i,j)\mid (i \in \b Z^{\geq 0})\land (j \in \nat)\} = \b d_2$. For each $m \in \b Z^{\geq 0}$ and $n \in \nat$, let
$H( m,n ) = \{\b f^2(i,j)\mid ( 0\leq i\leq m)\land (0\leq j \leq n)\land (i \in \b Z)\land(j \in \nat)\}.$ Then $H( m,n ) = d^2_{mn} \in {\cal F}'(\b d_2).$ Further, there is a $\b w^2_{mn} \in \b D_2,$ such that $H( m,n ) \subset \b S(\{\b w^2_{mn}\}).$ Hence, 
$$\forall y\forall x((y \in \b Z^{\geq 0})\land (x\in \nat) \to \exists z((z \in \b D_2)\land$$ 
$$(H( y,x ) \in {\cal F}'(\b d_2))\land(H( y,x ) \subset \b S(\{z\})))),\eqno (3)$$
holds in $\cal M$, hence, the *-transfer of (3) holds in $\hyper {\cal M}.$ This yields
$$\forall y\forall x((y \in \Hyper {\b Z}^{\geq 0})\land (x\in \hypernat) \to \exists z((z \in \Hyper {\b D_2})\land$$ 
$$(\hyper H( y,x ) \in \Hyper  {\cal F}'(\hyper {\b d_2}))\land(\hyper H( y,x ) \subset \Hyper {\b S}(\{z\})))).\eqno (4)$$
The conclusions follow as in Theorem 4.1 and this completes the proof. \qed\parm

\noindent{\bf Theorem 4.3.} {\it Consider primitive time interval $3 = (-\infty,0].$ Interval $3$ is partitioned into  intervals $[c_i, c_{i+1}),\ldots, [c_{-2},c_{-1}),[c_{-1},c_0], \ i \in \b Z^{\leq 0},\ i <-2.$ Let $\b d_3$
 be a developmental paradigm defined on $R_3$. For each infinite $\lambda \in \nat_\infty$ 
and $\mu \in \b Z_\infty^{\leq 0} = \Hyper {\b Z}^{\leq 0} - \b Z^{\leq  0},$ 
there is a $w^3_{\mu\lambda} \in \Hyper {\b D_3}$ 
and an hyperfinite $d^3_{\mu\lambda} \subset \hyper {\b d}_3$ 
such that $\b d_3 \subset d^3_{\mu \lambda}\subset \Hyper {\b S}(\{w^3_{\mu\lambda}\}).$ 
 Further, $x\in d^3_{\mu\lambda}$ if and only if there exist an $i \in \Hyper {\b Z}^{\leq 0},\  \mu \leq i$ and $j \in \hypernat, \ 0\leq j \leq \lambda,$ such that $x =\Hyper {\b f^3}(i,j).$} \parm\vfil\eject

\noindent{\bf Theorem 4.4.} {\it Consider primitive time interval $4 = (-\infty,+\infty).$ Interval $4$ is partitioned into  intervals $[c_{i},c_{i+1}),\ i \in \b Z$. Let $\b d_4$ be a developmental paradigm defined on $R_4$. 
For each $\lambda \in \nat_\infty$, $\nu \in \b Z_\infty^{\leq 0} = \Hyper {\b Z}^{\leq 0} - \b Z^{\leq 0}$ and $\gamma \in \b Z^{\geq 0}_\infty = \Hyper {\b Z}^{\geq 0} - \b Z^{\geq 0},$ there is a $w^4_{\nu \gamma\lambda} \in \Hyper {\b D^4}$ and an hyperfinite $d^4_{\nu \lambda \gamma} \subset \hyper {\b d_4}$ such that $\b d_4 \subset d^4_{\nu\gamma \lambda}\subset \Hyper {\b S}(\{w^4_{\nu \gamma \lambda}\}).$
Further, $x\in d^4_{\nu\gamma\lambda}$ if and only if there exist an $i \in \Hyper {\b Z},\ \nu \leq i \leq 0,$ and $k \in \Hyper {\b Z},\ 0\leq k \leq \gamma,$ and $j \in \hypernat,\  j\leq \lambda$ such that $x =\Hyper {\b f^4}(i,j)$ or $x =\Hyper {\b f^4}(k,j).$} \pars

Proof. The convergence requirement implies that $\{\b f^4(i,j)\mid (i \in \b Z^{\leq 0})\land(j \in \nat)\} \cup \{\b f^4(k,j)\mid (k \in \b Z^{\geq 0})\land(j \in \nat)\} = \b d_4.$ For $m \in \b Z^{\leq 0},\ p \in \b Z^{\geq 0}$ and $n \in \nat,$ let 
$H( m,p, n ) = \{\b f^4(i,j)\mid ( i\geq m)\land(i \in \b Z^{\leq 0})\land(j \leq n)\land (j \in \nat)\} \cup  \{\b f^4(k,j)\mid (k\leq p)\land (k \in \b Z^{\geq 0})\land(j \in \nat)\land(j \leq n)\}.$ Then $ H( m,p, n ) \in {\cal F}'(\b d_4).$ Further, there is a $\b w^4_{mpn} \in \b D_4$ such that $H( m,p,n ) \subset \b S(\{\b w^4_{mpn}\}).$ 
Expressing these conclusions formally yields the results  and this completes the proof.    \qed\parm

It is useful to investigate the actual objects contained in $\Hyper {\b S}(\{x\}),$ where $x$ is any of the ultrawords determined by the above theorems. This has been done for ultrawords generated in a slightly different manner and the result is stated in Theorem 10.1.1 in [7]. However, substitute $\b D_q$ for $\b M_{\b d_q} - \b d_q$  throughout that proof. Then, for each of the four types of intervals, the value of $\Hyper (H \cdots )$ for specific members of $\Hyper {\b Z_q} \times \hypernat$ is a hyperfinite developmental paradigm. The members of a developmental paradigm $\r d_q$ are considered as propositional atoms and three informal sets that correspond to the  $\hyper {\b A},\ Q_q ,\ d^q$ are shown to be disjoint in the altered proof for Theorem 10.1.1 [7]. This yields the following theorem. \parm 

{\bf Theorem 4.5} { \it For each $q  = 1,2,3,4$, let $w_q \in \Hyper {\b D_q}$ be an ultraword that exists by Theorem 2.q and let $\b d_q$ be the corresponding developmental paradigm and $d^q$ the corresponding hyperfinite set, where $d^q \subset \hyper {\b d_q}$ and $\b d_q \subset d^q .$ 
Then $\Hyper 
{\b S}(\{w_q\}) = \hyper {\b A} \cup Q_q \cup d^q,$ where for internal hyperfinite
$d^q,\ \b d_q \subset d^q \subset \Hyper {\b d_q}$ and 
internal $Q_q$ is composed of hyperfinite $(\geq 1)$ 
conjunctions $($i.e. $i({\rm \sp and\sp})$ $)$ of distinct members of $d_q$ and $w_q \in Q_q.$ Further$,$ each member of 
$d^q$ and no other member is used to form the hyperfinite 
conjunctions in $Q_q$ and members of $d^q$ are the only members of $\Hyper {\b S(\{w_q\})}$ 
without a special conjuction and $w_q$ is a hyperfinite conjuction, without repetition, of the members of $d^q$. Moreover, $\hyper {\b A},\ Q_q$ and $d^q$ are mutually disjoint.}\parm 
Each of the above five theorems is applicable to ``instructions or rules.'' For this case, in the W statement the word ``description'' is replaced with the phrase ``instruction'' or a similar term. The ultrawords that exist for the instructional rules are usually denoted by $\#w.$ \parm\vfil\eject

\noindent {\bf 5. Logic-System Signatures}. \par\medskip 
In formal logic, a certain amount of mental activity must be done before a formal proof is presented. For example, in most cases of interest, one needs to select finitely many well-formed formulas (wwfs) from potentially-infinite collections of wwfs. This is an acceptable process as modeled by a finite choice function. Further, such things as whether a variable is free or bound may need to be determined and when generalization is appropriate. Of course, there is also the mental activity required just to represent a collection of symbols in the proper form. When a formal deduction is presented, none of this mental activity is presented, although it might be discussed in an external manner using a metalanguage. Thus, not exhibiting such mental activity in the final product is a basic mathematical approach. In what follows, such external mental activity is also required and not represented in the final results. \par\smallskip

For a given nonempty language L, science-community scientific theories are discussed in [4]. One considers an implicit or explicit general rules of reference $\rm RI(L)$ that generates a finite consequence operator 
$\rm S_N$ that represents a particular scientific theory. For each $\rm X\subset L$, even with the realism relation $\rm R_N(X) = S_N(X) - X$ applied, there is a vast amount of extraneous ``deduction'' where the deduced members of 
$\rm R_N(X)$ are used to obtain the actual ``descriptions, words or images'' as a subset of $\rm R_N(X)$ that can be perceived. The term ``perceived'' often means ``to become aware of, through application of a set of defined human or machine sensory apparatus.''\par\smallskip

{\it The actual set $\rm P \subset L$  that constitutes what is termed here as ``perceived'' or ``observed'' should be explicitly defined by a science-community for a specific scientific theory or physical law.} Hence, for each $\rm X \subset P$,  $\rm P_X=P\cap S_N(X)= P_N(X)$ is the ``deduced'' perceived entities. Further, if a statistical statement is included that implies that members of $\rm P_X$ only have a certain probability of being perceived, then the ultralogic investigated in [5] is coupled with the $\rm P_X$ images. A ``signature'' is an entity that signifies the presence of a specific process or object. The operator $\rm P_N,$ defined on the set of all subsets of P, $\rm \power {P},$  is a finite consequence operator. \par\smallskip

The notion of the {\bf J}-relation as defined in [4] is now modified. Certain members of $\rm X$ may need to be tagged if they are also members of $\rm P_X.$ These members of $\rm X$ are considered as not being altered by the physical processes involved. The modified $\bf J'$-binary relation behaves like an identity relation for members of $\rm X$ except that the second coordinate is the same as the first coordinate with one additional fixed symbol attached to each member and the symbol does not appear in any of the perceived members of P. This symbol would not affect the actual ``meaning'' of any perceived member of L except that the symbol indicates that no change has been made in the expression denoted by the symbol by the physical processes being modeled. Note that, in what follows, physical laws are considered as producing a theory via a collection of rules of inference. \par\smallskip

For a given nonempty $\rm X \subset P$, a ``behavior-signature'' (${\rm B_X}$-signature) and the ``theory (or physical law)-signature'' (${\rm RI(P_N)}$-signature) are determined by ${\rm P_N}$. Note: In many of these investigations, the customary notation for ``n-tuples'' is employed where the actual definition may require the more formal definition by the ordered pair concept and induction or functions defined on various $[1,n],\ n >0, \ n\in \nat.$ \parm

{\bf Definition 5.1.} Given perceived $\rm P\subset L$ and a nonempty finite $\rm \{x_1,\ldots,x_n\} =X \subset P.$ If $\rm P_N(X) - X \not= \emptyset,$ define a {\bf behavior-signature} as $\rm {B_X}= \{(x_1,\ldots,x_n,x_{n+1})\mid x_{n+1} \in P_N(X) - X\}.$ Define the {\bf theory-signature} to be the unification $\rm {\rm RI(P_N)} = \bigcup\{R_X\mid (\emptyset \not= X \in {\cal F}(P)\},$ where $\cal F$ is the finite power set operator. \par\medskip

Definition 5.1 is equivalent to the rules of inference $\rm RI^*(P_N)$ as defined in [1, p. 204]. The difference is that the definition in [1] has the finite subsets of $\r P$, more or less, gathered together relative to cardinality and $\r R^1 = \emptyset$. Hence, Theorem 2.4 in [1] applies. Thus, the finite consequence operator generated by $\rm RI(P_N)$ is $\rm P_N.$ Although only objects different from members of X are used to obtain $\rm B_X$, $\rm X \subset P_N(X).$\parm 

The behavior-signatures are a more refined notion in that they are more specifically associated with deductive thought. 
When the realism relation is applied to $\rm P_N(X)$, then $\rm X$ is removed.  This does not remove the tagged members of $\rm P_N(X).$ The statements of the physical laws, logical axioms and other extraneous material not considered as members of $\rm P$ are now removed as ``deduced'' entities. Assume that $\rm S_N$ models a physical theory in that the theory processes are defined on all nonempty subsets of a language $\rm L_1 \subset L$  and the theory faithfully predicts the behavior of entities as they are described by members of $\rm L_1$. Let $\rm P\subset L_1.$ Hence, $\rm P_N$  affects specific perceived  objects $\rm X\subset P$ and yields perceived objects that may or may not differ in some describable sense from the original $\rm X\subset P$. The relation $\rm \{(X, P_N(X)-X)\mid  (X \subset P)\land (P_N(X)  \not= \emptyset)\}$ is contained in a ``physical process relation'' (See reference [8].)  Indeed, physical science-communities  attempt to show that it is equal to the physical process relation. \pars

There is a type of converse to Definition 5.1. Rather than starting with the $\rm S_N$, one can use observations and consider selecting a nonempty finite observation $\rm \{x_1,\ldots,x_n\} = X \subset P$. Assume physical processes applied to X yields a perceived $\rm X' \subset P.$  For each $\rm X,\ X'$, let $\rm B'_X = \{(x_1,\ldots,x_n,x_{n+1})\mid x_{n+1} \in X'\}.$   One considers  a ``unification'' $\rm \bigcup\{B'_X\mid (\emptyset \not= X \in{\cal F}(L)\}= {\rm RI'}.$ However, due to the general logic-system algorithm, even if one considers the finite logic-systems $\rm B'_X$ as separately applied, there are examples where the results need not be the same as those obtain by application of $\rm RI'.$  This fact can have significance for empirical science, where only such behavior-signatures are used to establish a rational theory $\rm S_N$. There are various reasons for this such as not knowing which objects in X are actually altered by the physical processes. One approach to correct this problem is to analyze carefully the data produced, alter how the data are expressed and produce a collection of  behavior-signatures that do correspond to those obtained from the corresponding $\rm RI'.$ In this case, the $\rm RI'$ can be consider as a representation for a physical law. Of course, these signature ideas may be applied to other appropriate ``natural'' laws that may not be considered as satisfying the strict definition for what constitutes a physical law. \par\medskip

\noindent {\bf 6. The Extended Standard Part Operator.}\par\medskip

One of the most significant operators used within nonstandard analysis is the ``standard part operator,'' ${\tt st}\colon G(0) \to \real,$ where $G(0)$ is the set of all ``finite'' (``limited'') hyperreal numbers [2, p. 17]. It is a point function that is defined in [1] via a set function extension. Although there are other procedures that lead to the standard real values, for subparticle representations of the form $(a_1,a_2,a_3,a_4, \ldots )= (k,\lambda,a_3,a_4,\ldots)$ [7, p. 99], the standard part is $\st {(k,\lambda,a_3,\cdots)} = (0,0,\st {a_3},\st {a_4}, \cdots)).$ [Note: In nonstandard analysis, there are two (isomorphic) ways to define n-tuples. One of these, that can be used for the formal definition for subparticle representation, where there are finitely many or denumerable many coordinates, is via the usual set of all functions $f$ from the indexing set $\emptyset \not= A\subset \nat$ to $\bigcup\{A_i\mid i \in A\}$ such that $f(i) \in A_i,$ where $A_1 = A_2 = \hypernat,\ A_i = G(0), i\geq 3.$] For a given application, let SP denote that set of all subparticle representations. \par\smallskip

For $\tt st$, let the corresponding set operator defined on each subset of SP be denoted by $\tt St.$ That is, for any $A \subset {\rm SP},\ \St {A} =\{\st {x} \mid x \in A\}.$  The operator $\tt st$ has the property that it is (composition) idempotent on members of $G(0).$ Let $X \subset \rm SP.$ If $X = \emptyset,$ then $\St {\emptyset} = \{\st {y}\mid y \in \emptyset\} = \emptyset= \St {\St {X}}.$ Suppose $X \not= \emptyset.$ Then $ \St {\St {X}}= \St {\{\st {x}\mid x \in X\}}= \{\st {\st {x}}\mid x \in X\}= \{\st {x}\mid x \in X\} = \St {X}.$ Hence, $\tt St$ is idempotent on $\power {{\rm SP}}.$ \par\smallskip

 For each $X \in \power {{\rm SP}},$ let $'\St {X} = X \cup \St {X}.$ 
The map $\tt St$ is closed under union. Hence, ${'\St {'\St {X}}}= {'\St {X \cup \St {X}}} = (X \cup \St {X}) \cup (\St {X \cup \St {X}})=  (X \cup \St {X}) \cup (\St {X} \cup 
\St {\St {X}}) = (X \cup \St {X})\cup (\St {X} \cup \St {X}) = {'\St {X}}.$ It follows that, for each $X \in \power {{\rm SP}},$ (2, 3) $X \subset {'\St {'{\St {X}}}}= {'\St {X}}\subset {\rm SP}.$ 
For $X \in \power {{\rm SP}},$ let ${\cal F}(X)$ denote the set of all finite subsets of $X$ and $x \in {'{\tt St(X)}} = X \cup \St {X}.$ Since $\tt st$ is a point map from which $\tt St$ is defined, then there exists $\{z_0\} \in \power {X}$ such that $x \in \{z_0\}\cup \St{\{z_0\}} \subset \bigcup \{F \cup \St {F}\mid F \in {\cal F}(X)\}= \bigcup\{'\St {F}\mid F \in {\cal F}(X)\}\subset X \cup \St {X}= {'{\tt St(X)}}.$ Thus (4) ${'{\tt St(X)}}= \bigcup\{{'\St {F}}\mid F \in {\cal F}(X)\}.$ Hence, for SP, $'\tt St$ satisfies axioms 2, 3, 4 in [7, p. 12] for the finite consequence operator. This operator is equivalent to a general logic-system [1]. Notice that for each $\emptyset \not= Y \in \power {{\rm SP}}$ such that for each $f \in Y$ and $i \geq 3,\ f(i)\in G(0) - \real$  the corresponding realism relation $R(Y) = {'\St {Y}} - Y = \St {Y}.$ As modeled by  *-linear transformations [7, p. 4, last paragraph], among other procedures, such $Y$ are precisely those used for a major step in GGU-model physical-entity generation. The operator ${'\tt ST}$ can also be used for the GID-model interpretation [3]. \par\medskip 

\centerline{\bf References}\par\medskip
\noindent [1] Herrmann, Robert A., (2006). General logic-systems and finite consequence operators, Logica Universalis 1:201-208. Portions appear in \hfil \break http://arxiv.org/abs/math/0512559 )\par\smallskip
\noindent [2] Herrmann, Robert A., (2003). Nonstandard Analysis - A Simplified Approach.\hfil\break http://arxiv.org/abs/math/0310351 \par\smallskip
\noindent [3] Herrmann, Robert A., (2002). Science Declares Our Universe IS Intelligently Designed, Xulon Press, Fairfax VA. \par\smallskip
\noindent [4] Herrmann, Robert A., (2001). Hyperfinite and standard unifications for physical theories, Internat. J. Math. and Math. Sci., 28(2):93-102.\hfil\break http://arxiv.org/abs/physics/0105012 \par\smallskip
\noindent [5] Herrmann, Robert A., (2001). Ultralogics and probability models, Internat. J. Math. and Math. Sci., 27(5):321-325. \hfil\break
http://arxiv.org/abs/quant-ph/0112037 \par\smallskip
\noindent [6] Herrmann, Robert. A., (1994). Solutions to the ``General Grand Unification Problem,'' and the Questions ``How Did Our Universe Come Into Being?" and ``Of What is Empty Space Composes?", 
http://arxiv.org/abs/astro-ph/9903110 \par\smallskip
\noindent [7] Herrmann, Robert. A., (1979 - ). The Theory of Ultralogics,  \hfil\break 
Part I at http://arxiv.org/abs/math/9903081 \hfil\break Part II at http://arxiv.org/abs/math/9903082 \par\smallskip
\noindent [8] Herrmann, R. A. Evidence, http://www.serve.com/herrmann/evidence.htm\par
\bigskip
\noindent {\it E-mail address:} rah@usna.edu herrmann@serve.com

\end